\documentclass[a4paper,9pt]{article}
\usepackage[latin1]{inputenc}
\usepackage{graphicx}
\usepackage{amsmath}
\usepackage{amsthm}
\usepackage{amssymb}
\usepackage{color}
\usepackage{mathrsfs}
\usepackage{fancyhdr}
\usepackage{setspace}
\usepackage[center]{subfigure}
\usepackage{float}
\usepackage[center,small,bf,ruled]{caption}
\usepackage{vmargin}

\usepackage[T1]{fontenc}
\usepackage{pbsi}
\usepackage{aurical}
\usepackage{calligra}
\usepackage{yfonts}
\allowdisplaybreaks

\setmarginsrb{25mm}{25mm}{25mm}{25mm}
             {0mm}{10mm}{0mm}{10mm}

\begin{document}

%%%%%%%%%%%%%%%%%%%%%%%%%%%%%%%%%%%%%%%%%%%%%%%%%%

\long\def\symbolfootnote[#1]#2{\begingroup
\def\thefootnote{\fnsymbol{footnote}}\footnote[#1]{#2}\endgroup} %per mettere le footnote con i simboli

\numberwithin{equation}{section}

\newtheorem {teo} {Theorem} [section]
\newtheorem {cor} {Corollary} [section]
\newtheorem {lem}[teo] {Lemma}

\theoremstyle{definition}
\newtheorem {oss} {Remark} [section]

\newcommand {\Dim} {\textsc{Proof\\ }}
\newcommand {\Fine} { $\blacksquare$ \\}
\newcommand{\q}{\left}
\newcommand{\p}{\right}
\date{} % delete this line to display the current date
\title{ {\bf Angular processes related to Cauchy\\ random walks }} 
\maketitle
\author{
{\center
{ \large V. Cammarota \symbolfootnote[2]{Dipartimento di Statistica, Probabilit\`a e Statistiche applicate, University of Rome `La Sapienza',  P.le Aldo Moro 5, 00185 Rome, Italy.  Tel.: +390649910499, fax: +39064959241. {\it E-mail address}: valentina.cammarota@uniroma1.it.} \;\;\;\;\;  E. Orsingher \symbolfootnote[3]{Corresponding author. Dipartimento di Statistica, Probabilit\`a e Statistiche applicate, University of Rome `La Sapienza',  P.le Aldo Moro 5, 00185 Rome, Italy.  Tel.: +390649910585, fax: +39064959241. {\it E-mail address}: enzo.orsingher@uniroma1.it. }} \\
\vspace{5mm}
}
\vspace{1cm}

\begin{abstract}
We study the angular process related to random walks in the Euclidean and in the non-Euclidean space where steps are Cauchy distributed. 

This leads to different types of non-linear transformations of Cauchy random variables which preserve the Cauchy density. We give the explicit form of these distributions for all combinations of the scale and the location parameters. 

Continued fractions involving Cauchy random variables are analyzed. It is shown that the $n$-stage random variables are still Cauchy distributed with parameters related to Fibonacci numbers. This permits us to show the convergence in distribution of the sequence to the golden ratio. 
\end{abstract}


AMS Classification 60K99\\

\section{Introduction}

We consider a particle starting from the origin $O$ of $\mathbb{R}^2$ which takes initially a horizontal step of length $1$ and a vertical one, say ${\mathbb C}_1$,   with a standard Cauchy distribution. It reaches therefore the position $(1,{\mathbb C}_1)$. The line $l_1$ joining the origin with $(1,{\mathbb C}_1)$ forms a random angle $\Theta_1$ with the horizontal axis (See Figure \ref{eucl}). 

On $l_1$ the traveller repeats the same movement with a step of unit length (either forward or backward) along $l_1$ and a standard Cauchy distributed step, say ${\mathbb C}_2$, on the line orthogonal to $l_1$. The right triangle obtained with the last two displacements has an hypothenuse belonging to the line $l_2$ with random inclination $\Theta_2$ on $l_1$.\\
\begin{figure}[h]
 \centering
    {\includegraphics[width=14.5cm, height=5.5cm]{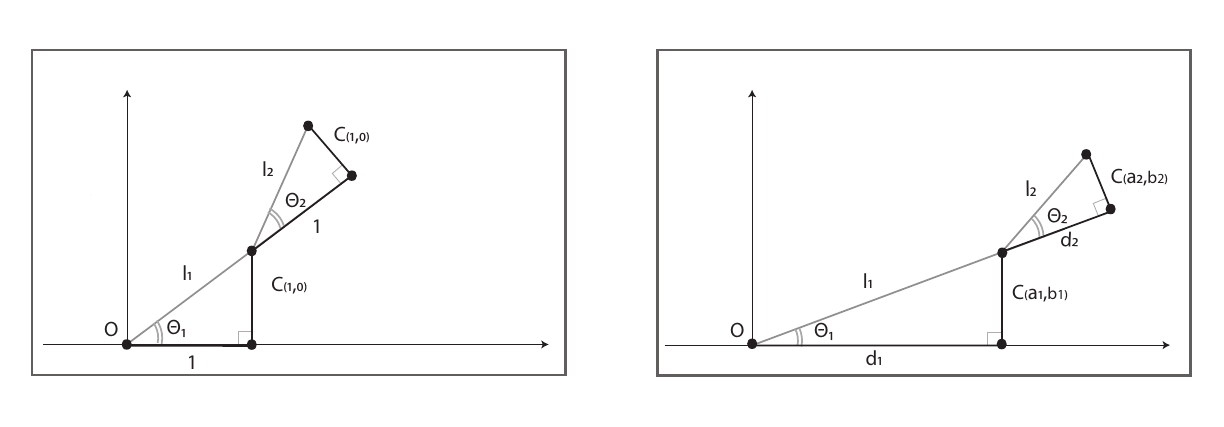}}
 \caption{The angular process in the Euclidean plane. By $C_{(a_j,b_j)}$ we indicate the j-th random displacement with Cauchy distribution possessing scale parameter $a_j$ and location parameter $b_j$.} \label{eucl}
 \end{figure}
After $n$ steps the sequence of random angles $\Theta_1,\cdots,\Theta_n$ describes the rotation of the moving particle around the starting point, their partial sums describe an angular random walk which can be written as 
 \begin{equation} \label{uno} {S}_{n}=\Theta_1+\cdots+\Theta_n=\sum_{j=1}^{n}\arctan {\mathbb C}_{j}\end{equation}
where ${\mathbb C}_{j}$ are independent standard Cauchy random variables. If the random steps of the planar random walk above were independent Cauchy random variables with scale parameter $a_j$ and location parameter $b_j$ then the process (\ref{uno}) must be a little bit modified and rewritten as 
\begin{equation} \label{due}
S_{n}=\Theta_1+\cdots+\Theta_n=\sum_{j=1}^{n}\arctan C_j,
\end{equation}
where $C_j \sim C_{(a_j,b_j)}$. The model (\ref{due}) can be extended also to the case where the first step has length $d_j$ and the second one is Cauchy distributed with scale parameter $a_j$ and position parameter $b_j$ (see Figure \ref{eucl}), then 
$$\tan \Theta_j=C_{ \left(\frac{a_j}{d_j},\frac{b_j}{d_j}\right)}.$$

 The same random walk can be generated if the two orthogonal steps, at each displacement, are represented by two independent Gaussian random variables $X_j$ and $Y_j$. In this case, for each right triangle, we can write
$$\tan \Theta_j=\frac{X_j}{Y_j}.$$   
If $X_j$ and $Y_j$ are two standard independent Gaussian random variables then $\tan \Theta_j=\frac{X_j}{Y_j}$ possesses standard Cauchy distribution and we get the model in (\ref{uno}). The model (\ref{due}) can be obtained by considering orthogonal Gaussian steps with different variances and in this case the scale parameter of the random variables $C_{j}$ is the ratio $a_j=\sigma_j^Y/ \sigma_j^X$.\\

The model (\ref{uno}) describing the angular random process has an hyperbolic counterpart. We consider a particle starting from the origin $O$ of the Poincar\'e half-plane $\mathbb{H}^2_+=\{(x,y): \; y>0\}$. At the $j$-th displacement, $j=1,2\dots$, the particle makes two steps of random hyperbolic length $\eta_{j}$ and $\hat{\eta}_{j}$ on two orthogonal geodesic lines. The $j$-th displacement leads to a right triangle $T_j$ with sides of length $\eta_{j}$ and $\hat{\eta}_{j}$ and random acute angles $\Theta_{j}$ and $\widehat{\Theta}_{j}$.  In each triangle $T_j$ the first step is taken on the extension of the hypotenuse of the triangle $T_{j-1}$ (see Figure \ref{hyperb}). From hyperbolic trigonometry (for basic results on hyperbolic geometry see, for example, Faber \cite{Faber}) we have that 
$$\sin  \Theta_{j}=\frac{\sinh \hat{\eta}_{j}}{\sqrt{\cosh^2 \eta_{j} \cosh^2 \hat{\eta}_{j}-1}},\;\;\;\;       \cos \Theta_{j}=\frac{\sinh \eta_{j} \cosh \hat{\eta}_{j}}{\sqrt{\cosh^2 \eta_{j} \cosh^2 \hat{\eta}_{j}-1}}.$$
From the above expressions we have that 
$$\tan  \Theta_{j}=\frac{\tanh \hat{\eta}_{j}}{\sinh \eta_{j}}.$$
If we take independent random hyperbolic displacements $\eta_j$ and $ \hat{\eta}_{j}$ such that the random variables $E_j=\frac{\tanh \hat{\eta}_{j}}{\sinh \eta_{j}}$ are standard Cauchy distributed then $\Theta_{j}=\arctan {\mathbb C}_{j}$. If the triangles $T_j$ were isosceles then 
$\tan \Theta_{j}=\frac{1}{\cosh \eta_{j}}$
and the angle $\Theta_{j} \in [-\frac{\pi}{4}, \frac{\pi}{4}]$ so that in this case the Cauchy distribution cannot be attributed to $\tan \Theta_{j}$.\\

\begin{figure}[h]
 \centering
    {\includegraphics[width=14.5cm, height=5.5cm]{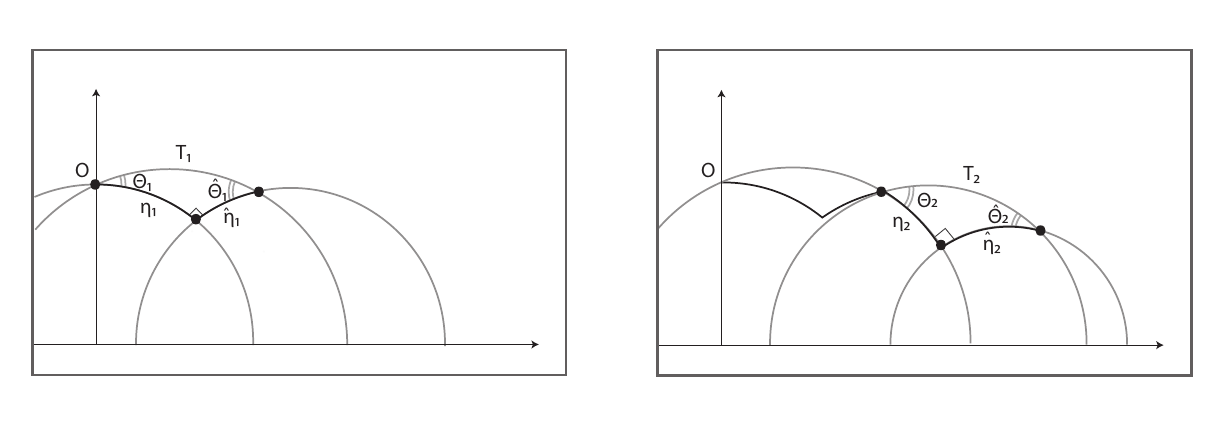}}
 \caption{The angular random process in the Poincar\'e half-plane.} \label{hyperb}
 \end{figure}

In the model described here the random steps (and therefore the random angular windings $\Theta_{j}$) are independent. If we consider the model of papers \cite{travel} and \cite{cascade}, where the displacements are taken orthogonally to the geodesic lines joining the origin $O$ of $\mathbb{H}^2_+$ with the positions occupied at deviation instants, the angular displacements $\Theta_{j}$ must be such that 
\begin{equation}
\sin \Theta_{j}=\frac{\sinh \eta_j }{\sqrt{1+\prod_{r=1}^{j} \cosh^2 \eta_r}}=\sinh \eta_j \cos \left(\arctan \prod_{r=1}^{j}\cosh ^2 \eta_r\right)  \nonumber
\end{equation} and therefore involve dependent random variables. 

For the area of the random hyperbolic triangle $T_j$ we note that
\begin{eqnarray*}
\mathrm{area}(T_j)&=&\frac{\pi}{2} - \Theta_{j}-\hat{\Theta}_{j}=\frac{\pi}{2}-\left[ \arctan \left( \frac{\tanh \hat{\eta}_{j}}{\sinh \eta_{j}} \right)+\arctan \left( \frac{\tanh \eta_{j}}{\sinh \hat{\eta}_{j}} \right) \right] \\
%&=&\frac{\pi}{2}-\arctan \left( \frac{\cosh \eta_{j}+\cosh \hat{\eta}_{j}}{\sinh \eta_j \sinh \hat{\eta}_j} \right)\\
&=&\frac{\pi}{2}-\arctan \left( \frac{\coth \eta_{j}}{\sinh  \hat{\eta}_{j}}+ \frac{\coth  \hat{\eta}_{j}}{ \sinh {\eta}_{j}} \right)=\mathrm{arcotan} \left( \frac{\coth {\eta}_{j} }{\sinh \hat{\eta}_{j}}+\frac{\coth \hat{\eta}_{j}}{\sinh {\eta}_{j}}  \right).
\end{eqnarray*}
Since each acute angle inside $T_j$ is linked to both sides of the triangle, the analysis of the random process $A_n=\sum_{j=1}^{n}\mathrm{area}(T_j)$ is much more complicated and we drop it.\\

Let $C_j\sim C_{(a_j,b_j)}$, $j=1,2\dots$ be independent Cauchy random variables where $a_j$ is the scale parameter and $b_j$ is the location parameter. In the study of the angular random walk (\ref{uno}) and (\ref{due}) we must analyze the distribution of the following non-linear transformations of Cauchy random variables: 
\begin{equation} \label{u}
U=\frac{C_1+C_2}{1-C_1C_2}
\end{equation} 
since 
\begin{equation*}
\arctan C_1+\arctan C_2= \arctan \frac{C_1+C_2}{1-C_1C_2}.
\end{equation*}

Since the Sixties a wide class of non-linear transformations of Cauchy random variables has been considered. Williams \cite{williams69}, Knight \cite{knight76} and  Letac \cite{letac77} proved that transformations of the form 
$$\epsilon f(x)=k x+\alpha-\sum_{i=1}^n \frac{\mu_i}{x-\gamma_i}$$
where $\epsilon=\pm1$, $\alpha, \gamma_i \in \mathbb{R}$ and $k, \mu_i \ge 0$ preserve the Cauchy distribution. 

In particular, in Williams \cite{williams69} it is proved the following characterization for Cauchy random variables. The random variable $X$ is a standard Cauchy if and only if $(1+bX)/(b-X)$ is a standard Cauchy for some constant $b$ which is not the tangent of a rational multiple of $\pi$. 

Knight \cite{knight76} asserts that a random variable $X$ is of Cauchy type if and only if the random variable $(aX+b)/(cX+d)$ is still of Cauchy type, whenever $ad-bc \ne 0$. 

Our problem is more strictly related to the results obtained by Pitman and Williams \cite{pitmanwilliams66}. They proved that the standard Cauchy distribution is preserved under certain types of transformations represented by meromorphic functions whose poles are all real and simple. As a corollary they obtained that, if $P$ and $Q$ are two independent random variables uniformly distributed in $(0,\pi)$, then the random variables $X=\tan P$ and $Y=\tan Q$ are standard Cauchy and 
$$\tan(P+Q)=\frac{X+Y}{1-XY}\stackrel{law}{=}X.$$

We will show that the random variable (\ref{u}) is endowed with Cauchy distribution in a much more general situation, namely when the random variables $C_j$, $j=1,2$, have non-zero location parameters $b_j$ and scale parameters $a_j$. The scale and location parameters of $U$ depend on both parameters $a_j$ and $b_j$ suitably combined.\\

In particular, if  $b_1=b_2=0$ and $a_1=a_2=1$, then $U$ is still distributed as a standard Cauchy distribution and therefore in (\ref{uno}) we have that $${S}_n\stackrel{law}{=}\arctan {\mathbb C}.$$
Since also $\frac{1}{{\mathbb C}}$ is a standard Cauchy (for basic properties of Cauchy random variables see, for example, Chaumont and Yor \cite{yor} page 105), from (\ref{u}), a number of related random variables preserving the form of the Cauchy distribution can be considered. For example, the following random variables 
$$Z_1=\frac{{\mathbb C}_1{\mathbb C}_2+1}{{\mathbb C}_1-{\mathbb C}_2},\;\;\;\;\;\;\;\;Z_2=\frac{1-{\mathbb C}_1{\mathbb C}_2}{{\mathbb C}_1+{\mathbb C}_2},\;\;\;\;\;\;\;\;Z_3=\frac{{\mathbb C}_1+{\mathbb C}_2}{{\mathbb C}_1{\mathbb C}_2-1},$$ 
also possess standard Cauchy distribution.\\
% We can also derive much more complicated random variables by suitably combining three (or more) independent standard Cauchy ${\mathbb C}_1, {\mathbb C}_2, {\mathbb C}_3$ random variables as 
%$$Z_4=\frac{{\mathbb C}_1+\frac{{\mathbb C}_2+{\mathbb C}_3}{1-{\mathbb C}_2 {\mathbb C}_3}}{1-{\mathbb C}_1\frac{{\mathbb C}_2+{\mathbb C}_3}{1-{\mathbb C}_2 {\mathbb C}_3}}=\frac{{\mathbb C}_1+{\mathbb C}_2+{\mathbb C}_3-{\mathbb C}_1{\mathbb C}_2 {\mathbb C}_3}{1-{\mathbb C}_1{\mathbb C}_2-{\mathbb C}_1 {\mathbb C}_3-{\mathbb C}_2{\mathbb C}_3}$$
%and so on.\\

If $b_1=b_2=0$ and the scale parameters $a_1, a_2$ are different, then (\ref{u}) still preserves the Cauchy distribution but with scale parameter equal to $\frac{a_1+a_2}{1+a_1a_2}$ and location parameter equal to zero.  This can be grasped by means of the following relationship 
\begin{equation} \label{ai}
\arctan C_{1}+\arctan C_{2}\stackrel{law}{=}\arctan \left\{\frac{a_1+a_2}{1+a_1a_2} {\mathbb C} \right\},
\end{equation}
where $C_j\sim C_{(a_j,0)}$. Result (\ref{ai}) is illustrated in Figure \ref{cann}. 

\begin{figure}[h]
 \centering
    {\includegraphics[width=12.5cm, height=6cm]{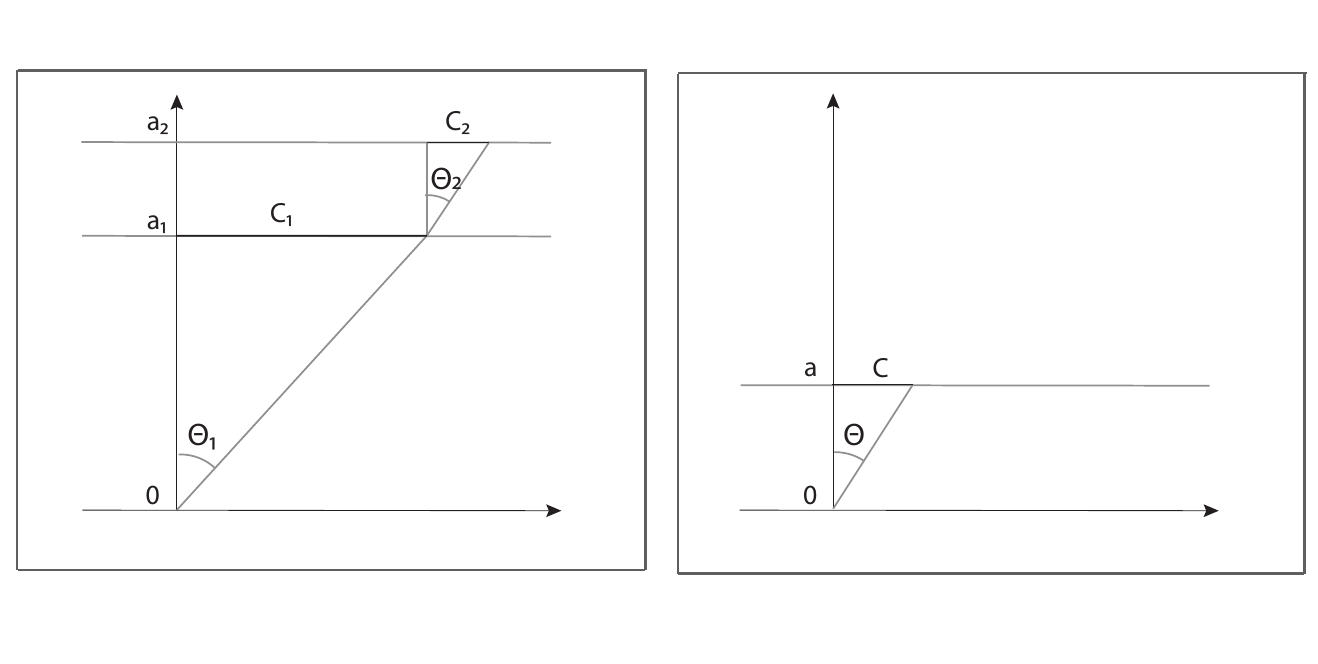}}
 \caption{The figure shows that shooting a ray with inclination $\Theta_1$, uniformly distributed, against the line at distance $a_1$ and then shooting a ray with a uniformly distributed angle $\Theta_2$ on the line at distance $a_2$ is equivalent to shooting on the barrier at the distance $a=\frac{a_1+a_2}{1+a_1a_2}$ with a uniformly distributed angle $\Theta$.} \label{cann}
 \end{figure}

By iterating the process (\ref{ai}) we arrive at the formula 
\begin{equation*}
\sum_{j=1}^{3} \arctan C_{j}\stackrel{law}{=}\arctan \left\{ \frac{\sum_{j=1}^{3} a_j +a_1a_2a_3 }{1+\sum_{i\ne j} a_ia_j} {\mathbb C} \right\}
\end{equation*}
which gives an insight into further extensions of the process outlined above.\\
%Many other relationships can be produced by combining the above results and we can observe that if $C_1\sim C_{(a_1,0)}$ and $C_2\sim C_{(a_2,0)}$ are independent Cauchy random variables, then 
%$$Z_5=a_1a_2 \frac{C_1+C_2}{C_1C_2-(a_1a_2)^2} $$
%also is a centered Cauchy random variable with scale parameter equal to $\frac{a_1+a_2}{1+a_1a_2}$.\\

Much more complicated are the cases where the location parameters of the Cauchy distributions are different from zero. For the special case where $C_1$ and $C_2$ are independent Cauchy such that $C_1\sim C_{(1,b)}$ and $C_2\sim C_{(1,b)}$, the random variable (\ref{u}) still possesses Cauchy density with scale parameter $\frac{2b^2+4}{b^4+4}$ and location parameter $\frac{2b^3}{b^4+4}$. \\ 

We have obtained the general distribution of (\ref{u}) where $C_1$ and $C_2$ are independent Cauchy such that $C_1\sim C_{(a_1,b_1)}$ and $C_2\sim C_{(a_2,b_2)}$ and also the distribution of 
$${U}=\frac{\gamma C_1+\delta C_2}{\alpha-\beta C_1C_2}$$
for arbitrary positive real numbers $\beta, \gamma, \delta$ and $\alpha+\beta \ne 0$. In particular, if $C_1$ and $C_2$ are independent standard Cauchy random variables then $U$ is Cauchy with scale parameter equal to $\left|\frac{\gamma+\delta}{\alpha+\beta}\right|$ and location parameter equal to zero.\\

In the last section we have examined continued fractions involving Cauchy random variables. In particular we have studied 
\begin{equation} \label{ln}
V_n=\frac{1}{1+\frac{1}{1+\frac{1}{1+\cdots \frac{1}{1+{\mathbb C}}}}}
\end{equation}
and 
\begin{equation} \label{gn}
U_n=\frac{1}{1+\frac{1}{1+\frac{1}{1+\cdots \frac{1}{1+{\mathbb C}^2}}}}
\end{equation}
which generalize the random variables 
$V_1=\frac{1}{1+{\mathbb C}}$ and $U_1=\frac{1}{1+{\mathbb C}^2}$. Continued fractions involving random variables have been analyzed in Chamayou and Letac \cite{cham} and more recently in Asci, Letac and Piccioni \cite{picc}. The random variable $U_1$ has the arcsine distribution in $[0,1]$, while $U_t=tU_1$, with $t>0$, has distribution 
$$\mathrm{Pr}\{U_t \in ds\}=\frac{\mathrm{d}s}{\pi\sqrt{s(t-s)}}, \hspace{1cm}0<s<t.$$\\
For each $n\ge1$, the random variables $V_n$, are Cauchy distributed with scale parameter $a_n$ and  position parameter $b_n$ that can be expressed in terms of Fibonacci numbers. This permits us to prove the monotonicity of $a_n$ and $b_n$ and that $\lim_{n \to \infty}a_n = 0$ and $\lim_{n \to \infty }b_n = \phi-1$ where $\phi=\frac{1+\sqrt{5}}{2}$ is the golden ratio. Finally we obtain that the sequence of random variables $1+V_n$ and $1+U_n$, $n \ge 1$, converges in distribution to the number $\phi=\frac{1+\sqrt{5}}{2}$. This should be expected since it has the infinite fractional expansion
\begin{equation}
\frac{1+\sqrt{5}}{2}=1+\frac{1}{1+\frac{1}{1+\cdots}}
\end{equation}
 which is related to (\ref{ln}) and (\ref{gn}).

\section{Centered Cauchy  random variables}

In this section we study the distribution of the following random variable 
\begin{equation}  \label{UUU} U=\frac{\gamma{\mathbb C}_1+\delta{\mathbb C}_2}{\alpha-\beta {\mathbb C}_1{\mathbb C}_2},\end{equation}
where ${\mathbb C}_1$ and ${\mathbb C}_2$ are independent, standard Cauchy. We assume, without restrictions, that  $\beta, \gamma, \delta$ are non-negative real numbers all different from zero (because of the symmetry of ${\mathbb C}_j$, $j=1,2$) and  $\alpha+\beta \ne 0$. In the next theorem we prove that $U$ is still Cauchy distributed.

%The first approach is based on the evaluation of some cumbersome integrals while the second one is hinged on well-known properties of the Cauchy distribution. In the second proof of Theorem \ref{one} it is relevant to note that  \begin{equation} \label{pitman} V=\frac{{\mathbb C}_1+{\mathbb C}_2}{1- {\mathbb C}_1{\mathbb C}_2},\end{equation}
%is a standard Cauchy. This is proved in classical papers like Pitman and Williams \cite{pitmanwilliams66} and can be shown also by direct calculations.

\begin{teo} \label{one}
The random variable $U$ in (\ref{UUU}) possesses Cauchy distribution with scale parameter equal to $\left|\frac{\gamma+\delta}{\alpha+\beta}\right|$ and position parameter equal to zero. We can also restate the result in symbols as 
\begin{equation*}\label{UU} U\stackrel{law}{=}\frac{\gamma+\delta}{\alpha+\beta} {\mathbb C}.\end{equation*}
\end{teo}
\textsc{Proof 1}\\
The density of the random variable $U$ can be obtained by means of the transformation 
\begin{equation*}
\left\{
\begin{array}{lr} u=\frac{\gamma x+\delta y}{\alpha - \beta x y},\\
v=y,
\end{array}
\right.
\end{equation*}
with $(x,y) \in \mathbb{R}^2$ and Jacobian equal to 
$$|J|=\frac{\alpha \gamma+\delta \beta v^2}{(\gamma+\beta u v)^2}.$$
The joint density $g=g(u,v)$ reads 
\begin{eqnarray*}g(u,v)&=& \frac{1}{\pi^2(1+v^2)} \frac{1}{1+ \left( \frac{\alpha u- \delta v}{\gamma+\beta u v}  \right)^2} \frac{ \alpha \gamma+\delta \beta v^2}{(\gamma+\beta u v)^2} \hspace{1cm} (u,v)\in \mathbb{R}^2\\
&=&  \frac{1}{\pi^2(1+v^2)} \frac{\alpha \gamma+\delta \beta v^2}{v^2(\beta^2 u^2+\delta^2)+2 u v (\beta \gamma - \delta \alpha)+\alpha^2 u^2+\gamma^2 }
\end{eqnarray*}
and must be integrated with respect to $v$ in order to obtain the distribution of $U$. Therefore

%Since 
%\begin{eqnarray*}
 %\mathrm{Pr}\left\{\frac{\gamma {\mathbb C}_1+ \delta {\mathbb C}_2}{\alpha-\beta {\mathbb C}_1 {\mathbb C}_2}  < u \right \}&=& \mathrm{E}\left\{  \mathrm{Pr} \left\{ \left. \frac{\gamma {\mathbb C}_1+ \delta {\mathbb C}_2}{\alpha-\beta {\mathbb C}_1 {\mathbb C}_2}<u \right| {\mathbb C}_2  \right\}     \right\}=\frac{1}{\pi^2} \int_{-\infty}^{\infty} \frac{\mathrm{d}v}{1+v^2}  \int_{-\infty}^{\frac{\alpha u- \delta v }{\gamma+ \beta u v}} \frac{\mathrm{d}x  }{1+x^2}, 
 %\end{eqnarray*}
%we arrive at the following integrals 
\begin{eqnarray} \label{2.9}
 \mathrm{Pr}\left\{\frac{\gamma {\mathbb C}_1+ \delta {\mathbb C}_2}{\alpha-\beta {\mathbb C}_1 {\mathbb C}_2} \in \mathrm{d} u \right \}%&=&\left\{\frac{\mathrm{d}}{\mathrm{d} u} \frac{1}{\pi^2} \int_{-\infty}^{\infty} \frac{\mathrm{d}v}{1+v^2}  \int_{-\infty}^{\frac{\alpha u- \delta v }{\gamma+ \beta v u}} \frac{\mathrm{d}x  }{1+x^2} \right\} \mathrm{d} u\nonumber \\
&=&\frac{\mathrm{d} u}{\pi^2} \int_{-\infty}^{\infty}\q[ \frac{1}{1+v^2} \frac{\alpha \gamma+\delta \beta v^2}{v^2(u^2 \beta^2+\delta^2)+2vu(\beta \gamma-\delta \alpha)+u^2\alpha^2+\gamma^2}\p]\mathrm{d}v\\
&=&\frac{\mathrm{d} u}{\pi^2} \int_{-\infty}^{\infty}\q[ \frac{Av+B}{1+v^2}+\frac{Cv+D}{v^2(u^2 \beta^2+\delta^2)+2vu(\beta \gamma-\delta \alpha)+u^2\alpha^2+\gamma^2} \p]\mathrm{d}v  \nonumber
\end{eqnarray}
where
\begin{equation*}
\left\{
\begin{array}{lr} 
\vspace{0.3cm}
A=\frac{2u(\beta \gamma-\alpha \delta)(\beta \delta-\alpha \gamma)}{[u^2(\alpha-\beta)^2+(\gamma-\delta)^2][u^2(\alpha+\beta)^2+(\gamma+\delta)^2]},\\
\vspace{0.3cm}
B=\frac{(\alpha \gamma-\beta \delta)[(\gamma^2-\delta^2)+u^2(\alpha^2-\beta^2)]}{[u^2(\alpha-\beta)^2+(\gamma-\delta)^2][u^2(\alpha+\beta)^2+(\gamma+\delta)^2]},\\
\vspace{0.3cm}
C=-A(u^2\beta^2+\delta^2),\\
\vspace{0.3cm}
D=\frac{(\beta \gamma-\alpha \delta)[u^4\alpha \beta(\beta^2-\alpha^2 )+u^2(\beta \gamma-\alpha \delta)(3 \alpha\gamma-\beta \delta)+\gamma \delta(\gamma^2-\delta^2)]}{[u^2(\alpha-\beta)^2+(\gamma-\delta)^2][u^2(\alpha+\beta)^2+(\gamma+\delta)^2]}.
\end{array}
\right.
\end{equation*}
We start by evaluating the first part of the integral (\ref{2.9}):
\begin{eqnarray} \label{2.10}
&& \int_{-\infty}^{\infty} \left[  \frac{Av}{1+v^2}+  \frac{Cv}{v^2(u^2 \beta^2+\delta^2)+2vu(\beta \gamma-\delta \alpha)+u^2\alpha^2+\gamma^2}\right] \mathrm{d}v \nonumber\\
&=&\frac{A}{2}\int_{-\infty}^{\infty} \left[ \frac{2v}{1+v^2} - \frac{2v (u^2\beta^2+\delta^2 ) \pm 2u(\beta\gamma-\delta\alpha)}{v^2(u^2 \beta^2+\delta^2)+2vu(\beta \gamma-\delta \alpha)+u^2\alpha^2+\gamma^2}\right] \mathrm{d}v \nonumber\\
&=&\frac{A}{2}\left. \lim_{d \to \infty, c\to -\infty} \log \left(\frac{1+v^2}{v^2(u^2 \beta^2+\delta^2)+2vu(\beta \gamma-\delta \alpha)+u^2\alpha^2+\gamma^2}\right) \right|_{c}^{d}\nonumber \\
&&+Au (\beta\gamma-\delta\alpha)\int_{-\infty}^{\infty}\frac{1}{v^2(u^2 \beta^2+\delta^2)+2vu(\beta \gamma-\delta \alpha)+u^2\alpha^2+\gamma^2} \mathrm{d}v\nonumber \\
&=&Au (\beta\gamma-\delta\alpha)\int_{-\infty}^{\infty}\frac{1}{v^2(u^2 \beta^2+\delta^2)+2vu(\beta \gamma-\delta \alpha)+u^2\alpha^2+\gamma^2} \mathrm{d}v \nonumber \\
&=&Au (\beta\gamma-\delta\alpha) \frac{\pi}{u^2\alpha\beta+\gamma\delta},
\end{eqnarray}
where the last integral is obtained by means of the change of variable
$$v \sqrt{u^2\beta^2+\delta^2}+\frac{u(\beta\gamma-\delta\alpha)}{\sqrt{u^2\beta^2+\delta^2}}=z\sqrt{u^2 \alpha^2+\gamma^2-\frac{u^2(\beta\gamma-\delta\alpha)^2}{u^2\beta^2+\delta^2}}=z \frac{u^2 \alpha \beta+\gamma \delta}{\sqrt{u^2 \beta^2+\delta^2}}.$$
In view of result (\ref{2.10}) and inserting the values of $A$, $B$ and $D$ we have that 
\begin{eqnarray*}
&&\frac{\mathrm{d} u}{\pi^2} \int_{-\infty}^{\infty}\q[ \frac{Av+B}{1+v^2}+\frac{Cv+D}{v^2(u^2 \beta^2+\delta^2)+2vu(\beta \gamma-\delta \alpha)+u^2\alpha^2+\gamma^2} \p]\mathrm{d}v\\
&=&\frac{\mathrm{d}u}{\pi}\frac{1}{u^2\alpha\beta+\gamma\delta}[Au(\beta\gamma-\alpha\delta)+B(u^2\alpha\beta+\gamma\delta)+D]\\
&=&\frac{\mathrm{d}u}{\pi}\left[ \frac{(\gamma \delta+u^2 \alpha \beta)(\beta \gamma-\alpha \delta)[u^2(\beta^2-\alpha^2)+(\gamma^2-\delta^2)]   }{(u^2\alpha\beta+\gamma\delta)[u^2(\alpha-\beta)^2+(\gamma-\delta)^2][u^2(\alpha+\beta)^2+(\gamma+\delta)^2]}  \right.\\
&&+\left. \frac{(\gamma \delta+u^2 \alpha \beta)(\alpha \gamma-\beta \delta)[u^2(\alpha^2-\beta^2)+(\gamma^2-\delta^2)]   }{(u^2\alpha\beta+\gamma\delta)[u^2(\alpha-\beta)^2+(\gamma-\delta)^2][u^2(\alpha+\beta)^2+(\gamma+\delta)^2]}  \right]\\
&=&\frac{\mathrm{d}u}{\pi} \frac{u^2(\beta^2-\alpha^2)(\beta-\alpha)(\gamma+\delta)+(\gamma^2-\delta^2)(\alpha+\beta)(\gamma-\delta)  }{[u^2(\alpha-\beta)^2+(\gamma-\delta)^2][u^2(\alpha+\beta)^2+(\gamma+\delta)^2]}  \\
&=&\frac{\mathrm{d}u}{\pi} \frac{(\alpha+\beta)(\gamma+\delta)[u^2(\alpha-\beta)^2+(\gamma-\delta)^2] }{[u^2(\alpha-\beta)^2+(\gamma-\delta)^2][u^2(\alpha+\beta)^2+(\gamma+\delta)^2]}  \\
&=&\frac{\mathrm{d}u}{\pi} \frac{(\alpha+\beta) (\gamma+\delta)}{u^2(\alpha+\beta)^2+(\gamma+\delta)^2}.
\end{eqnarray*}

Another approach is based on the conditional characteristic function,  
\begin{eqnarray*}
\mathrm{E}\{e^{i \beta U}|\mathbb{C}_2=v\}&=&\frac{1}{\pi}\int_{-\infty}^\infty e^{i \beta (\frac{\gamma u+\delta v}{\alpha-\beta u v})} \frac{\mathrm{d}u}{1+u^2}\\
&=&\frac{1}{\pi} \int_{-\infty}^\infty e^{i \beta w} \frac{\alpha \gamma+\delta \beta v^2}{v^2(\beta^2 w^2+\delta^2)+2 v w (\beta \gamma-\delta \alpha)+w^2 \alpha^2+\gamma^2} \mathrm{d}w
\end{eqnarray*}
where $w=\frac{\gamma u+\delta v}{\alpha-\beta u v}$. The inverse Fourier transform gives the conditional density $g(w|v)$ and thus we arrive again at the integral (\ref{2.9}). \\

\Fine

\begin{oss}
A special case implied by Theorem \ref{one} concerns the random variable 
$$V=\frac{C_1+C_2}{1-C_1 C_2}$$
where $C_j\sim C_{(a_j,0)}$, with $j=1,2$, are independent Cauchy random variables with location parameter equal to zero and scale parameter $a_j>0$. If we choose $\gamma=a_1$, $\delta=a_2$, $\alpha=1$ and $\beta=a_1 a_2$, we conclude that 
$$V\stackrel{law}{=} \frac{a_1 {\mathbb C}_1+a_2 {\mathbb C}_1}{1-a_1 a_2 {\mathbb C}_1 {\mathbb C}_2}\stackrel{law}{=}\frac{a_1+a_2}{1+a_1 a_2}{\mathbb C}.$$
In the same way we obtain that 
\begin{eqnarray*}
\widehat{V}=\frac{\gamma C_1+ \delta C_2}{\alpha-\beta C_1C_2}\stackrel{law}{=} \frac{\gamma a_1 {\mathbb C}_1+\delta a_2 {\mathbb C}_2}{\alpha-\beta a_1a_2 {\mathbb C}_1{\mathbb C}_2}\stackrel{law}{=} \frac{\gamma a_1+\delta a_2}{\alpha+\beta a_1a_2} {\mathbb C}.
\end{eqnarray*}  
where $\alpha+\beta a_1 a_2 \ne 0$.
\end{oss}

 \begin{oss}
In view of Theorem \ref{one}  we can obtain  by recurrence the distribution of the current angle after $n$ steps for the angular random walk 
 $$S_n=\sum_{j=1}^n \arctan C_j$$
described in the introduction. We have 
$$S_n=\sum_{j=1}^{n}\arctan C_{j(a_j,0)}=\arctan C_{(\hat{a}_{n-1},0)}+\arctan C_n(a_{n},0)$$
 where $\arctan C(\hat{a}_{n-1},0)$ is the random variable $S_{n-1}$. In particular, if $a_j=1$ for $j=1,\dots n$, we have the following property of the standard Cauchy random variables ${\mathbb C}_j$
 \begin{equation*}
 \sum_{j=1}^{n} \arctan {\mathbb C}_{j}\stackrel{law}{=}\arctan {\mathbb C}.
 \end{equation*}
 \end{oss}

  \begin{oss}
 A simple byproduct of Theorem \ref{one} is that 
 \begin{eqnarray} \label{inte}
 E e^{i \beta U}&=&\frac{1}{\pi^2} \int_{\mathbb{R}^2} e^{i \beta \frac{x+y}{1-xy}} \frac{a_1a_2 \mathrm{d}x \mathrm{d}y}{(a_1^2+x^2)(a_2^2+y^2)} =\frac{1}{\pi^2} \int_{-\frac{\pi}{2}}^{\frac{\pi}{2}} \int_{-\frac{\pi}{2}}^{\frac{\pi}{2}} e^{i \beta \frac{a_1 \tan \theta_1+a_2\tan \theta_2}{1-a_1 a_2 \tan \theta_1 \tan \theta_2}} \mathrm{d}\theta_1 \mathrm{d}\theta_2 \nonumber \\
 &=&e^{-\frac{a_1+a_2}{1+a_1a_2}|\beta|}.
\end{eqnarray} 
In (\ref{inte}) we have used the transformations $x=a_1 \tan\theta_1$  and $y=a_2 \tan\theta_2$. In the special case $a_1=a_2=1$ the relationship (\ref{inte}) yields 
 \begin{eqnarray}\label{byp} 
 e^{- |\beta|}&=&\frac{1}{\pi^2} \int_{-\frac{\pi}{2}}^{\frac{\pi}{2}} \int_{-\frac{\pi}{2}}^{\frac{\pi}{2}} e^{i \beta \frac{\tan \theta_1+\tan \theta_2}{1-\tan \theta_1 \tan \theta_2}} \mathrm{d}\theta_1 \mathrm{d}  \theta_2 =\frac{1}{\pi^2} \int_{-\frac{\pi}{2}}^{\frac{\pi}{2}} \int_{-\frac{\pi}{2}}^{\frac{\pi}{2}} e^{i \beta \tan(\theta_1+\theta_2)} \mathrm{d}\theta_1 \mathrm{d}\theta_2\nonumber \\
 &=&\frac{2}{\pi^2}\int_{0}^{\pi} x \cos(\beta \tan x ) \mathrm{d}x.
\end{eqnarray} 
In the last step of (\ref{byp}) we have used the transformations $\theta_1+\theta_2=x$ and $\theta_2=y$. The integral (\ref{byp}) shows that, if $(\Theta_1, \Theta_2)$ is uniform in the square $S=\{(\theta_1, \theta_2): -\frac{\pi}{2}<|\theta_i|<\frac{\pi}{2},\; i=1,2\}$, then the random variable $W=\tan(\Theta_1+\Theta_2)$ has characteristic function $e^{-|\beta|}$ because $\Theta_1+\Theta_2$ is uniform and therefore $W$ is Cauchy distributed.
 \end{oss}
 
 \begin{oss}
It is well-known that  for a planar Brownian motion $\{(B_1(t), B_2(t)), t>0\}$ starting from $(x,y)$  the random variable $B_1(T_y)$ possesses Cauchy distribution with parameters $(x,y)$ where 
 $$T_y=\inf \{t>0:\; B_2(t)=0\}.$$
If the starting points of two planar Brownian motions $(B^i_1(t), B^i_2(t))$, for $i=1,2$, are located on the $y$ axis as in the Figure \ref{brow} then we have that 
\begin{eqnarray*}
\Theta=\Theta_1+\Theta_2&=&\arctan B^1_1(T_{a_1})+\arctan B^2_1(T_{a_2})\\ &\stackrel{law}{=}&\arctan C_1 +\arctan C_2\\&\stackrel{law}{=}&\arctan \frac{a_1+a_2}{1+a_1a_2} {\mathbb C}.
\end{eqnarray*}
where $C_1$ and $C_2$ are two independent Cauchy random variables with scale parameters $a_1$ and $a_2$ respectively and position parameter equal to zero. Therefore if the starting point of a third Brownian motion has coordinates $\left( 0, \frac{a_1+a_2}{1+a_1a_2} \right)$ then $B\left(T_{ \frac{a_1+a_2}{1+a_1a_2}} \right)$  represents its hitting position on the $x$-axis. This point forms with $(0,1)$ and the origin a right triangle with an angle equal to $\Theta=\Theta_1+\Theta_2$. 
\begin{figure}[h]
 \centering
    {\includegraphics[width=14.5cm, height=5.5cm]{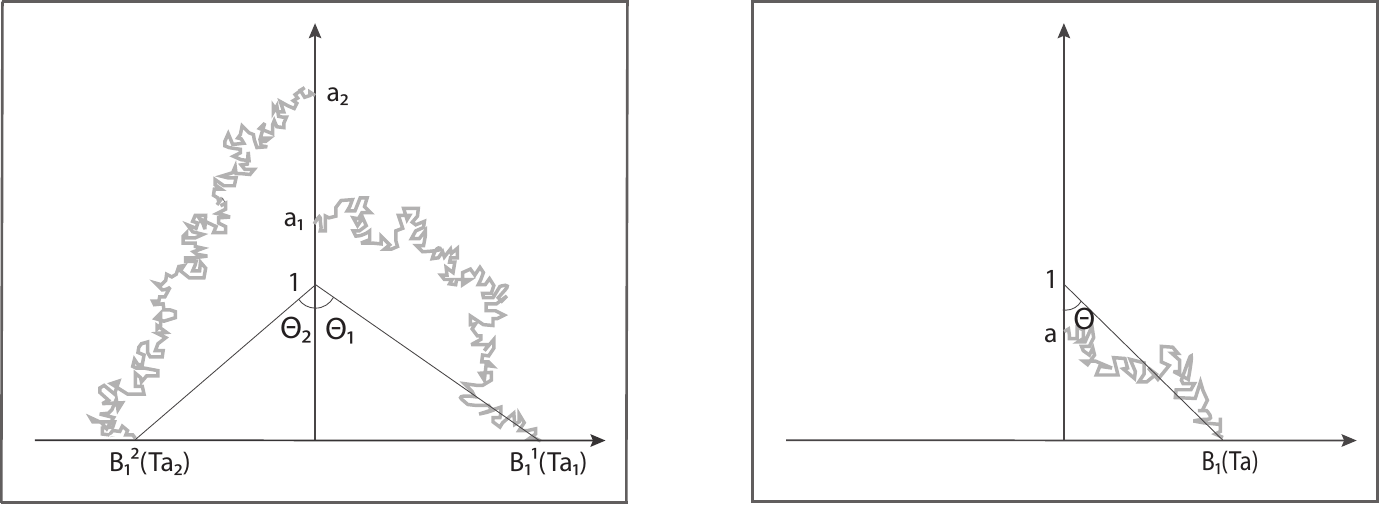}}
 \caption{The hitting position on the $x$-axis of a planar Brownian motion is Cauchy distributed. In the figure the random angles $\Theta_1$, $\Theta_2$ and $\Theta=\Theta_1+\Theta_2$ are shown. For $a=\frac{a_1+a_2}{1+a_1a_2}$ the right-hand figure shows the hitting position $B_1(T_a)$.} \label{brow}
 \end{figure}
\end{oss}

\section{Non-Centered Cauchy random variables}

For independent Cauchy random variables $C_1$ and $C_2$, with location parameters $b_1$ and $b_2$ and scale parameters $a_1$ and $a_2$, the random variable $U$ in (\ref{UUU}) is still Cauchy distributed with both parameters affected by the values of the location parameters $b_1$, $b_2$ and the scale parameters $a_1$, $a_2$.

\begin{teo} \label{g}
If $C_{i}$, $i=1,2$, are two independent, Cauchy random variables with location parameters $b_i$ and scale parameters $a_i$, then the random variable $U$ is still Cauchy distributed with scale parameter $$a_U=\frac{ (a_1+a_2)(1+a_1 a_2-b_1b_2)+(b_1+b_2)(a_1b_2+a_2 b_1)}{(1+a_1a_2-b_1b_2)^2+(a_1 b_2+a_2 b_1)^2 }, $$  and position parameter $$b_U=\frac{(a_1+a_2)(a_1b_2+a_2 b_1)-(b_1+b_2)(1+a_1 a_2-b_1b_2)}{(1+a_1a_2-b_1b_2)^2+(a_1 b_2+a_2 b_1)^2}.$$ 
 \end{teo}
 \Dim
 We obtain the density function of the random variable $U$ by observing that 
 \begin{equation*}
 \mathrm{Pr} \{U \in \mathrm{d}w\}= \mathrm{E}\left\{ \mathrm{Pr}\left\{\left. \frac{C_{1}+C_{2}}{1-C_{1} C_{2}} \in \mathrm{d}w \right| C_2\right\} \right\}
 \end{equation*} 
 and remarking that 
 \begin{equation*}
 \mathrm{Pr}\left\{\left. \frac{C_{1}+C_{2}}{1-C_{1} C_{2}} \in \mathrm{d}w \right| C_2=y\right\}=\frac{a_1 }{\pi \left[ \frac{w-y}{1+wy}-b_1   \right]^2+a_1^2} \frac{1+y^2}{(1+ w y)^2}  \mathrm{d}w.
 \end{equation*}
 Therefore
 \begin{eqnarray}
 \mathrm{Pr} \{U \in \mathrm{d}w \}&=&\frac{a_1 a_2 \; \mathrm{d}w}{\pi^2} \int_{-\infty}^{\infty} \frac{1 }{ \left[ \frac{w-y}{1+wy}-b_1   \right]^2+a_1^2} \frac{1+y^2}{(1+ w y)^2} \frac{1}{(y-b_2)^2+a_2^2}\mathrm{d}y \nonumber \\
 &=&\frac{a_1 a_2 \; \mathrm{d}w}{\pi^2} \int_{-\infty}^{\infty} \frac{1+y^2 }{ \left[ w-y-b_1 (1+ w y)   \right]^2+a_1^2 (1+ w y)^2} \frac{1}{(y-b_2)^2+a_2^2}\mathrm{d}y \nonumber \\
 &=& \frac{a_1 a_2 \; \mathrm{d}w}{\pi^2} \int_{-\infty}^{\infty} \frac{1+y^2 }{[y^2 E -2 y F +G] [y^2-2y H+K] } \mathrm{d}y, \label{jkjk}
 \end{eqnarray}
 where 
 \begin{eqnarray*}
 E=1+2 b_1 w+ (a_1^2+b_1^2)w^2, \hspace{1cm} F=-b_1-(a_1^2+b_1^2-1)w+b_1 w^2,
 \end{eqnarray*}
 \begin{eqnarray*}
 G=a_1^2+b_1^2-2 b_1 w +w^2,  \hspace{1cm} H=b_2,    \hspace{1cm} K=a_2^2+b_2^2.
 \end{eqnarray*}
 We rewrite the integral in (\ref{jkjk}) in the following form 
 \begin{eqnarray} \label{dsds}
  \mathrm{Pr} \{U \in \mathrm{d}w \}&=&\frac{a_1 a_2 \; \mathrm{d}w}{\pi^2} \int_{-\infty}^{\infty} \left[  \frac{A y+B}{y^2 E -2 y F +G} + \frac{ C y+D}{y^2 -2 y H +K} \right] \mathrm{d}y
 \end{eqnarray}
 with 
 \begin{equation} \label{fgfg}
\left\{
\begin{array}{lr} 
A+CE=0,\\
-2 A H+B-2 CF +DE=1,\\
AK-2 HB+CG-2 DF=0,\\
BK+DG=1.
\end{array}
\right.
\end{equation}
The integrals in (\ref{dsds}) can be worked out by means of the change of variables 
$$y \sqrt{E} - \frac{F}{\sqrt{E}}=x \sqrt{G-\frac{F^2}{E}}, \hspace{1cm} y-H=x \sqrt{K-H^2}.$$
The first integral becomes 
\begin{eqnarray*}
 \int_{-\infty}^{\infty}  \frac{A y+B}{ \left(y \sqrt{E} - \frac{F}{\sqrt{E}} \right)^2+G-\frac{F^2}{E}} \mathrm{d}y=\int_{-\infty}^{\infty} \frac{1}{E} \frac{BE+AF+x A \sqrt{GE-F^2}}{(x^2+1) \sqrt{GE-F^2}} \mathrm{d}x,
\end{eqnarray*}
and the second one takes the form  
\begin{eqnarray*}
 \int_{-\infty}^{\infty}  \frac{C y+D}{ \left(y - H \right)^2+K-H^2} \mathrm{d}y=\int_{-\infty}^{\infty} \frac{D+CH+x C \sqrt{K-H^2}}{(x^2+1) \sqrt{K-H^2}} \mathrm{d}x.
\end{eqnarray*}
A substantial simplification can be obtained because 
$$\sqrt{GE-F^2}=a_1 (1+w^2), \hspace{1cm} \sqrt{H-K^2}=a_2.$$
If we turn back to the distribution (\ref{dsds}), in view of the above calculations, we have that 
\begin{eqnarray*} \mathrm{Pr} \{U \in \mathrm{d}w \}&=&\frac{a_2 \mathrm{d}w}{\pi^2 (1+w^2) E}  \int_{-\infty}^{\infty} \frac{BE+AF+x A a_1 (1+w^2)}{x^2+1} \mathrm{d}x \\&&+ \frac{a_1 \mathrm{d}w}{\pi^2} \int_{-\infty}^{\infty} \frac{D+CH+x C a_2 }{1+x^2} \mathrm{d}x.
\end{eqnarray*}
We observe that, in view of the first equation of (\ref{fgfg}), we have that 
$$\frac{ a_1 a_2\; \mathrm{d}w}{\pi^2}  \left( \frac{A}{E}+C  \right)   \int_{-\infty}^{\infty} \frac{x}{1+x^2} \mathrm{d}x=0$$
and 
\begin{eqnarray*}
\mathrm{Pr} \{U \in \mathrm{d}w \}&=&\left\{ a_2 \frac{BE+AF}{\pi E (1+w^2)} +a_1 \frac{D+CH}{\pi}  \right\} \mathrm{d}w\\
&=& \left\{  a_2 \frac{B-CF}{\pi (1+w^2)} +a_1 \frac{D+CH}{\pi}  \right\} \mathrm{d}w.
\end{eqnarray*}
The values of $B$, $C$ and $D$ can be derived by solving the system (\ref{fgfg}), simplified as  
 \begin{equation*}
\left\{
\begin{array}{lr} 
B+2 C (HE-F)+DE=1,\\
-2 B H+C(G-EK)-2DF=0,\\
BK+DG=1.
\end{array}
\right.
\end{equation*}
By means of  cumbersome computations we arrive at the final result  
 \begin{eqnarray*}
&&\mathrm{Pr} \{U \in \mathrm{d}w \} \\
&&=\frac{1}{\pi} \frac{ (1+a_1 a_2-b_1b_2) [w(a_1 b_2+a_2 b_1)+a_1+a_2]- (a_1 b_2+a_2 b_1) [w(1+a_1 a_2-b_1b_2) -  (b_1+b_2)] }{[w(1+a_1 a_2-b_1b_2)-(b_1+b_2)]^2+[w(a_1 b_2+a_2 b_1)+a_1+a_2]^2}\\
&&= \frac{\frac{ (a_1+a_2)(1+a_1 a_2-b_1b_2)+(b_1+b_2)(a_1b_2+a_2 b_1)}{(1+a_1a_2-b_1b_2)^2+(a_1 b_2+a_2 b_1)^2 } }{\q[w+\frac{(a_1+a_2)(a_1b_2+a_2 b_1)-(b_1+b_2)(1+a_1 a_2-b_1b_2)}{(1+a_1a_2-b_1b_2)^2+(a_1 b_2+a_2 b_1)^2}\p]^2  + \q[  \frac{ (a_1+a_2)(1+a_1 a_2-b_1b_2)+(b_1+b_2)(a_1b_2+a_2 b_1)}{(1+a_1a_2-b_1 b_2)^2+(a_1 b_2+a_2 b_1)^2 }\p]^2}.
\end{eqnarray*}
%In the above calculations one must bear in mind that 
%\begin{eqnarray*}
%&&\frac{(a_1+a_2)^2+(b_1+b_2)^2}{(1+a_1a_2-b_1 b_2)^2+(a_1 b_2+a_2 b_1)^2 }- \q[  \frac{ (b_1+b_2)(a_1b_2+a_2 b_1)+(a_1+a_2)(1+a_1 a_2-b_1 b_2)}{(1+a_1a_2-b_1 b_2)^2+(a_1 b_2+a_2 b_1)^2 }\p]^2\\
%&&=\q[ \frac{(a_1b_2+a_2 b_1)(a_1+a_2)-(b_1+b_2)(1+a_1 a_2-b_1 b_2)}{(1+a_1a_2-b_1 b_2)^2+(a_1 b_2+a_2 b_1)^2} \p]^2
%\end{eqnarray*}
%because 
%\begin{eqnarray*}
%&&[(a_1+a_2)^2+(b_1+b_2)^2][(1+a_1a_2-b_1 b_2)^2+(a_1 b_2+a_2 b_1)^2]-[ (b_1+b_2)(a_1b_2+a_2 b_1)\\&&+(a_1+a_2)(1+a_1 a_2-b_1 b_2)]^2\\
%&=& (a_1+a_2)^2 (1+a_1a_2-b_1 b_2)^2+(b_1+b_2)^2(1+a_1a_2-b_1 b_2)^2+(a_1+a_2)^2 (a_1 b_2+a_2 b_1)^2\\&&+(b_1+b_2)^2(a_1 b_2+a_2 b_1)^2 -(b_1+b_2)^2(a_1b_2+a_2 b_1)^2-(a_1+a_2)^2 (1+a_1 a_2-b_1 b_2)^2\\&&-2(a_1+a_2)(b_1+b_2)(a_1b_2+a_2 b)(1+a_1 a_2-b_1 b_2)\\
%&=&[(a_1b_2+a_2 b_1)(a_1+a_2)-(b_1+b_2)(1+a_1 a_2-b_1 b_2)]^2 .
%\end{eqnarray*}
\Fine

\begin{oss}
In view of Theorem \ref{g} it is possible to obtain the following particular cases. 
\begin{itemize}
\item For $a_1=a_2=1$ and $b_1=b_2=b$, we have that 
$$a_U=\frac{2b^2+4}{b^4+4}, \hspace{2cm} b_U=\frac{2b^3}{b^4+4}.$$
\end{itemize}
This shows that $U$ has center of symmetry on the positive half-line if $b>0$ and on the negative half-line if $b<0$, therefore the non linear transformation $U$ preserves the sign of the mode.
\begin{itemize}
\item For $a_1=a_2=a$ and $b_1=b_2=b$ we have that 
$$a_U=\frac{2a(1+a^2+b^2)}{(1+a^2-b^2)^2+(2ab)^2}, \hspace{1cm}b_U=\frac{2b(a^2+b^2-1)}{(1+a^2-b^2)^2+(2ab)^2}.$$
\end{itemize}
We note that $a_U$ and $b_U$ depend simultaneously from the scale and location parameters of the random variables involved in $U$. 
\end{oss}

\section{Continued Fractions}

The property that the reciprocal of a Cauchy random variable has still a Cauchy distribution has a number of possible extensions which we deal with in this section.

We start by considering the sequence  
\begin{equation}  \label{vn}
V_1=\frac{1}{1+{\mathbb C}}, \hspace{1cm} V_2=\frac{1}{1+\frac{1}{1+{\mathbb C}}},\dots \hspace{1cm} V_n=\frac{1}{1+\frac{1}{1+\dots\frac{1}{1+{\mathbb C}}}},
\end{equation} 
and show the following theorem.

\begin{teo} \label{t4.1}
The random variables defined in (\ref{vn}) have Cauchy distribution $V_n\sim C_{(a_n,b_n)}$ where the scale parameters $a_n$ and the location parameters $b_n$ satisfy the recursive relationships
\begin{equation} \label{an} a_{n+1}=\frac{a_n}{(1+b_n)^2+a_n^2}, \hspace{2cm} n=1,2,\dots
\end{equation}
and
\begin{equation} \label{bn}b_{n+1}=\frac{b_n+1}{(1+b_n)^2+a_n^2}, \hspace{2cm} n=1,2,\dots\end{equation}
\end{teo}

\Dim 
Let us assume that $V_n$ possesses Cauchy density with parameters $a_n$ and $b_n$, therefore $V_{n+1}$ writes 
$$\mathrm{Pr}\{ V_{n+1}<v\}= \mathrm{Pr}\left\{ \frac{1}{1+V_n}<v\right\}=\mathrm{Pr} \left\{\frac{1}{1+a_n+b_n {\mathbb C}}<v\right\}.$$
After some computations the density of $V_{n+1}$ can be written as  
$$f_{V_{n+1}}(v)=\frac{\frac{a_n}{(1+b_n)^2+a_n^2}}{\pi \left[ v-\frac{b_n+1}{(1+b_n)^2+a_n^2}\right]^2+\frac{a_n^2}{[(1+b_n)^2+a_n^2]^2}}, \hspace{1cm} v \in \mathbb{R}.$$
It can be directly ascertained that $V_1$ possesses Cauchy distribution with parameters $a_1=1/2$ and $b_1=1/2$. 
\Fine

\begin{oss} \label{rem}
We have evaluated the following table of parameters $a_n$ and $b_n$:
$$\hspace{0.2cm}n\hspace{1.5cm}1\hspace{1.5cm}2\hspace{1.5cm}3\hspace{1.5cm}\dots \hspace{1.5cm} 10^{2}$$
$$\hspace{0.3cm}a_n\hspace{1.1cm}1/2\hspace{1.2cm}1/5\hspace{1.2cm}1/13\hspace{1cm}\dots\hspace{1cm}5.77 e^{-42}$$ 
$$\hspace{0.3cm}b_n\hspace{1.1cm}1/2\hspace{1.2cm}3/5\hspace{1.2cm}8/13\hspace{1cm}\dots\hspace{1cm}0.618034$$\\
For $n=1,2,3$ we can observe that the scale parameters $a_n$ coincide with the inverse of the odd-indexed Fibonacci numbers while the sequence $b_n$ has the numerators coinciding with the even-indexed Fibonacci numbers and the denominators correspond to the odd-indexed Fibonacci numbers.

In light of the previous considerations we can show that for $n\ge1$
\begin{equation} \label{Fn} b_n=\frac{F_{2n}}{F_{2n+1}},\hspace{2cm}a_n=\frac{1}{F_{2n+1}},\end{equation}
where $F_n$, $n\ge 0$ is the Fibonacci sequence. Recalling that the Fibonacci numbers admit the following representation (it can be easily checked by induction) 
\begin{equation} \label{fn}F_n=\frac{\phi^n-(1-\phi)^n}{\sqrt{5}}\end{equation}
where $\phi=\frac{1+\sqrt{5}}{2}$ is the golden ratio, we now prove that if $a_n$ and $b_n$ have the representation (\ref{Fn}), then also $a_{n+1}$ and $b_{n+1}$ can be expressed in the same form. From (\ref{an}) and (\ref{bn}) we have 
\begin{eqnarray*}
b_{n+1}&=&\frac{\frac{F_{2n}}{F_{2n+1}} +1}{\left(\frac{F_{2n}}{F_{2n+1}}+1 \right)^2 +\frac{1}{F^2_{2n+1}}}=\frac{F_{2n+2}F_{2n+1}}{F^2_{2n+2}+1}\\
&=&F_{2n+2} \left[ \frac{\phi^{2n+1}-(1-\phi)^{2n+1} }{\phi^{4n+4}+(1-\phi)^{4n+4}-2\phi^{2n+2}(1-\phi)^{2n+2}+5} \right]\sqrt{5}\\
&=&F_{2n+2} \left[ \frac{\phi^{2n+1}-(1-\phi)^{2n+1} } {\left[ \phi^{2n+1}-(1-\phi)^{2n+1}  \right]  \left[ \phi^{2n+3}-(1-\phi)^{2n+3}  \right]} \right] \sqrt{5}\\
&=&\frac{F_{2n+2}}{F_{2n+3}}.
\end{eqnarray*}
Similar calculations prove that $a_{n+1}=\frac{1}{F_{2n+3}}$. In view of representation (\ref{Fn}) and (\ref{fn}), it is easy to show that 
\begin{equation*}
\lim_{n \to \infty} b_n=\lim_{n \to \infty} \frac{F_{2n}}{F_{2n+1}}=\lim_{n \to \infty} \frac{1-\left(\frac{1-\phi}{\phi}  \right)^{2n}}{\phi-(1-\phi) \left(\frac{1-\phi}{\phi}  \right)^{2n} }=\frac{1}{\phi} =\phi-1, \hspace{1cm} \lim_{n \to \infty} a_n=\lim_{n \to \infty} \frac{1}{F_{2n+1}}=0.
\end{equation*} \\
Otherwise, observing that the sequence $b_n$, $n \ge 1$ is increasing, because 
\begin{eqnarray*}
\frac{b_{n+1}}{b_n}=\frac{F_{2n+2} F_{2n+1}}{F_{2n+3} F_{2n}}=\frac{\phi^{4n+3}+(1-\phi)^{4n+3}+1}{\phi^{4n+3}+(1-\phi)^{4n+3}-4}\ge 1
\end{eqnarray*}
and taking the limits in (\ref{an}) and (\ref{bn}) we have that 
\begin{equation} \label{lim}L=\frac{L}{(1+H)^2+L^2}, \hspace{2cm}H=\frac{H+1}{(1+H)^2+L^2},\end{equation}\\
where $H=\lim_{n \to \infty}a_n$ and $L=\lim_{n \to \infty}b_n$. From the relationships in (\ref{lim}) we derive the equality \begin{equation*} \label{lim2} \frac{L}{H}=\frac{L}{H+1}\end{equation*} that implies $L=0$. In fact, for $L\ne 0$, we arrive at the absurd that $H=H+1$. Substituting $L=0$ in the second formula of (\ref{lim}) we obtain $$H=\frac{1}{H+1},$$ since $H$ satisfies the algebraic equation $H^2+H-1=0$ it follows that $H=\phi-1$ where $\phi$ is the golden ratio (see Figure 5). \end{oss}

\begin{oss}
A slightly more general case concerns the sequence 
$$W_1=\frac{1}{c_1+d_1 C_{(a_0,b_0)}}=\frac{1}{c_1+a_0d_1+b_0d_1 {\mathbb C}},\hspace{0.5cm}W_2=\frac{1}{c_2+d_2 W_1},\hspace{0.5cm}W_3=\frac{1}{c_3+d_3 W_2},\dots$$
By performing calculations similar to those of Theorem \ref{t4.1} we have that $W_1$ has Cauchy distribution with scale parameter $a_1$ and position parameter $b_1$ such that $$a_1=\frac{d_1a_0}{(c_1+d_1b_0)^2+d_1^2a_0^2}, \hspace{0.5cm} b_1=\frac{c_1+d_1b_0}{(c_1+d_1b_0)^2+d_1^2a_0^2}.$$
Similarly, if $W_{n} \sim C_{(a_n,b_n)}$, than $W_{n+1}\sim C_{(a_{n+1},b_{n+1})}$ where
\begin{equation} \label{4.4} a_{n+1}=\frac{d_{n+1} a_n}{(c_{n+1}+d_{n+1}b_n)^2+d_{n+1}^2a_n^2}, \hspace{0.5cm} b_{n+1}=\frac{c_{n+1}+d_{n+1} b_n}{(c_{n+1}+d_{n+1}b_n)^2+d_{n+1}^2a_n^2}\end{equation}
for every $n \ge 2$. The sequences in (\ref{4.4}) for $c_n=d_n=1$ coincide with (\ref{an}) and (\ref{bn}).
\end{oss}

\begin{figure} \label{auchy^2}
 \centering
 %\subfigure[]
   {\includegraphics[width=7cm, height=7cm]{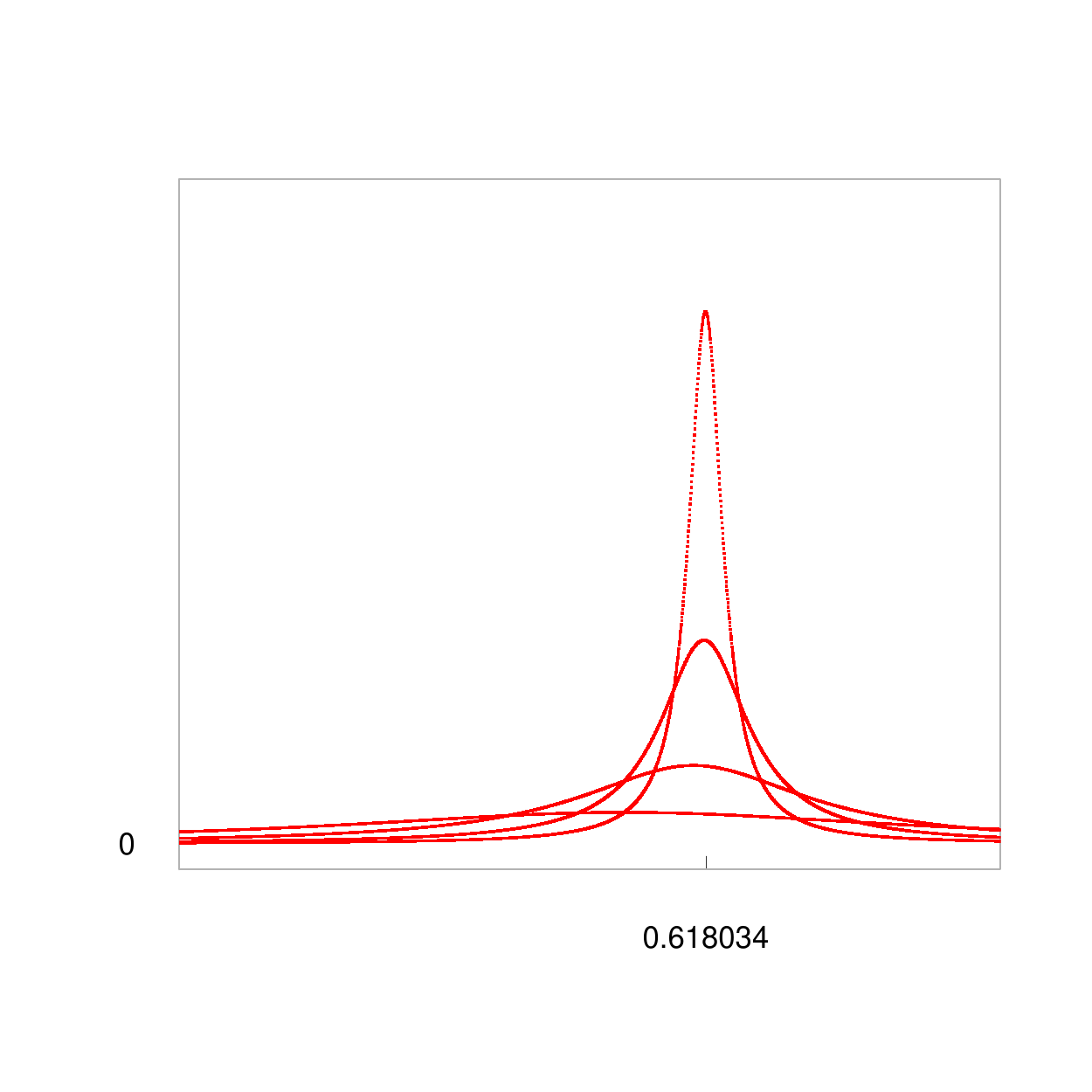}}
 \hspace{5mm}
 %\subfigure[]
   {\includegraphics[width=7cm, height=7cm]{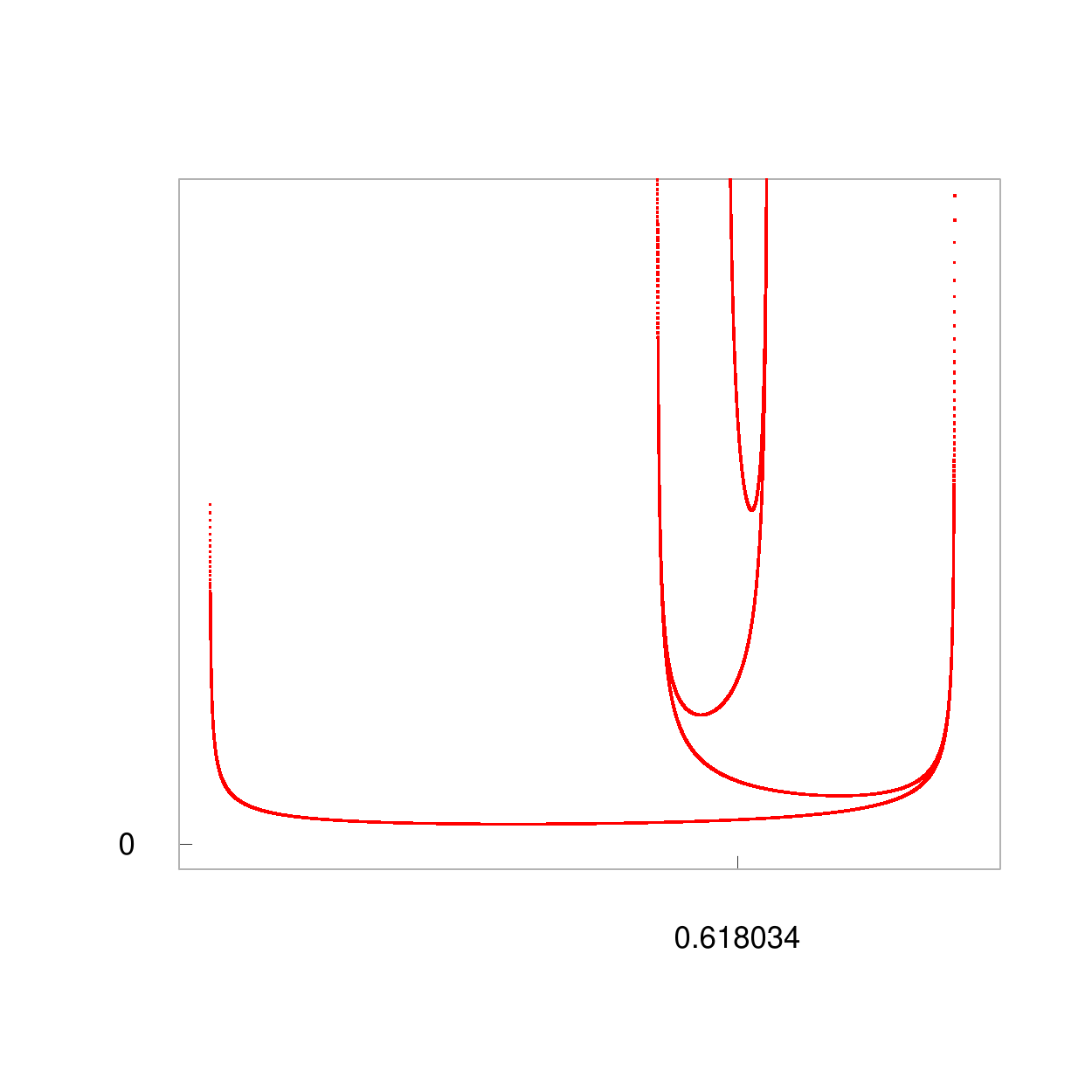}}
\caption{In the left figure the densities of the Cauchy random variables $V_1$, $V_2$, $V_3$ and $V_4$ are shown. In the right figure the densities of $U_1$, $U_2$, $U_3$ and $U_4$ are plotted.}
 \end{figure}

Another sequence of continued fractions involving the Cauchy distribution is the following one 
\begin{equation}  \label{un}
U_1=\frac{1}{1+{\mathbb C}^2}, \hspace{1cm} U_2=\frac{1}{1+\frac{1}{1+{\mathbb C}^2}},\dots \hspace{1cm} U_n=\frac{1}{1+\frac{1}{1+\dots\frac{1}{1+{\mathbb C}^2}}}
\end{equation} 
It is well-known that the random variable $U_1$ possesses the arcsin law.  Unlike the sequence $V_n$ studied above the sequence $U_n$, $n\ge 1$, has a density structure changing with $n$. Some calculations are sufficient to show that $U_1$, $U_2$, $U_3$, $U_4$ have density, respectively equal to 
\begin{eqnarray*}
f_{U_1}(u)&=&\frac{1}{\pi \sqrt{u(1-u)}}, \hspace{3.5cm}0< u <1,\\
f_{U_2}(u)&=&\frac{1}{\pi u \sqrt{(1-u)(2u-1)}}, \hspace{2.2cm}\frac{1}{2}< u <1,\\
f_{U_3}(u)&=&\frac{1}{\pi (1-u) \sqrt{(2u-1)(2-3u)}}, \hspace{1.2cm}\frac{1}{2}< u <\frac{2}{3},\\
f_{U_4}(u)&=&\frac{1}{\pi (2u-1) \sqrt{(2-3u)(5u-3)}}, \hspace{1cm}\frac{3}{5}< u <\frac{2}{3}.
\end{eqnarray*}
The general result concerning $U_n$ is stated in the next theorem.

\begin{teo}
For every $n \ge1$ the distribution of the random variable $U_{n}$ is given by
\begin{eqnarray} \label{dist}
\mathrm{Pr}\{U_n \in \mathrm{d}u\}&=&\frac{1}{\pi [(-1)^{n+1} \alpha_n+(-1)^{n}\beta_n u] } \frac{ 1}{\sqrt{(-1)^{n} \beta_n + (-1)^{n+1}(\alpha_n +\beta_n)u }} \nonumber \\
&& \times \frac{ 1}{\sqrt{(-1)^{n+1} (\alpha_n+\beta_n) + (-1)^{n}(\alpha_n +2\beta_n)u}} \mathrm{d}u,
\end{eqnarray}
where
$$(-1)^n \frac{\alpha_n + \beta_n}{\alpha_n+2 \beta_n} < (-1)^n u < (-1)^n \frac{\beta_n}{\alpha_n+\beta_n},$$
and  $\alpha_n, \beta_n \in \mathbb{N}$ satisfy the recursive relationships $\alpha_n=\beta_{n-1}$, $\beta_n=\alpha_{n-1}+\beta_{n-1}$.
\end{teo} 

\Dim
From (\ref{un}) we have that
$$U_{n+1}=\frac{1}{1+U_n},$$
then proceeding by induction, i.e. assuming that $U_n$ has distribution (\ref{dist}), we obtain that 
\begin{eqnarray*}
\mathrm{Pr} \{U_{n+1} \in \mathrm{d}u\}&=&\frac{\mathrm{d}}{\mathrm{d}u} \mathrm{Pr}\left\{U_n > \frac{1-u}{u}\right\} \mathrm{d}u= \frac{\mathrm{d}}{\mathrm{d}u} \int_{\frac{1-u}{u}}^{h(\alpha_n, \beta_n)}  \mathrm{Pr}\left\{U_n \in \mathrm{d}u\right\}\\
&=& \frac{1}{\pi} \frac{1}{u^2} \frac{1}{ [(-1)^{n+1} \alpha_n+(-1)^{n}\beta_n (\frac{1-u}{u})] } \frac{ 1}{\sqrt{(-1)^{n} \beta_n + (-1)^{n+1}(\alpha_n +\beta_n)(\frac{1-u}{u}) }} \nonumber \\
&& \times \frac{ 1}{\sqrt{(-1)^{n+1} (\alpha_n+\beta_n) + (-1)^{n}(\alpha_n +2\beta_n)(\frac{1-u}{u})}} \mathrm{d}u \\
&=&\frac{1}{\pi [(-1)^{n} \beta_n+(-1)^{n+1}(\alpha_n+\beta_n) u] } \frac{ 1}{\sqrt{(-1)^{n+1} (\alpha_n+\beta_n) + (-1)^{n}(\alpha_n +2\beta_n)u }} \nonumber \\
&& \times \frac{ 1}{\sqrt{(-1)^{n} (\alpha_n+2\beta_n) + (-1)^{n+1}(2\alpha_n +3\beta_n)u}} \mathrm{d}u.
\end{eqnarray*}
In the first integral the function $h(\alpha_n, \beta_n)$ represents the right boundary of the support of $U_n$. We conclude that $U_{n+1}$ possesses distribution (\ref{dist}) by taking $\alpha_{n+1}=\beta_n$ and $\beta_{n+1}=\alpha_n+\beta_n$. 
\Fine

\begin{oss}
The sequence $\beta_n $ is a Fibonacci sequence since we have that $\beta_n=\beta_{n-1}+\alpha_{n-1}=\beta_{n-1}+\beta_{n-2}$. We note that the sequence of coefficients $\alpha_n$ and $\beta_n$ are such that 
$$\lim_{n \to \infty}\frac{\alpha_{n+1}}{\beta_{n+1}}=\lim_{n \to \infty}\frac{\beta_{n}}{\beta_{n+1}} =\phi-1.$$
On the base of arguments similar to those of Remark \ref{rem} it is possible to show that the sequence $U_n$, $n \ge 1$, converges in distribution to $\phi-1$. In this case the upper and lower bounds of the domain of definition of the densities $f_{U_n}(u)$, $n \ge1$ are expressed as ratios of Fibonacci numbers (see Figure 5). 
%By means of observations similar to that  displayed in Remark \ref{rem} it is possible to show the convergence in distribution of the sequence of random variables $U_n$, $n \ge 1$, to the value $\phi-1$. In this case we get the sequence of ratios of Fibonacci numbers as upper and lower bounds of the domain of definition of the densities $f_{U_n}(u)$, $n \ge1$ (see Figure 5(b)).
\end{oss}
\vspace{1cm}
{\bf Acknowledgment} We are very grateful to the referee for his scholar report and also to have drawn our attention to some relevant references.

\end{document}